\documentstyle[leqno, amssymb, amsmath,12pt]{article}
\numberwithin{equation}{section}
\newtheorem{thm}{Theorem}[section]
\newtheorem{prop}[thm]{Proposition}
\newtheorem{lem}[thm]{Lemma}

\begin{document}

\centerline{\LARGE An analogue to the Witt identity \vspace{5mm}}
\centerline{\large G. A. T. F. da 
Costa \footnote{g.costa@ufsc.br}
and
G. A. Zimmermann \footnote{graciele@ifsc.edu.br}}
\centerline{\large Departamento de Matem\'{a}tica}
\centerline{\large Universidade Federal de Santa Catarina}
\centerline{\large 88040-900-Florian\'{o}polis-SC-Brasil}

\begin{abstract}

\vspace{5mm}

{\bf In this paper we 
solve
combinatorial and algebraic problems associated with a multivariate  identity first considered by S. Sherman which he called
{\it an analog to the Witt identity}. We extend previous results obtained for the univariate case.}

\end{abstract}


\vspace{10mm}

\section{Introduction} \label{sec : Int}

In \cite{ssherman} S. Sherman considered
the formal identity 
in the indeterminates 
$z_{1}, ...,z_{n}$:
\begin{equation} \label{A}
\prod_{m_{1},...,m_{R} \geq 0 }
(1 + z_{1}^{m_{1}}...z_{R}^{m_{R}})^{N_{+}}
(1 - z_{1}^{m_{1}}...z_{R}^{m_{R}})^{N_{-}}
=
\prod_{j=1}^{R} (1 + z_{j})^{2}
\end{equation}
where
$N_{+}$ and $N_{-}$
are the numbers of distinct classes of equivalence of nonperiodic closed paths
with positive and negative signs, respectively, which
traverse without backtracking 
$m_{i}$ times edge $i$, $i=1, ...,R$, of a graph $G_{R}$ 
with $R > 1$ edges forming loops counterclockwisely 
oriented and
hooked to a single  vertex, $\sum m_{i} \geq 1$. 

In \cite{ssherman} Sherman refers to equation \eqref{A} as 
{\it an analog to the Witt identity}. The reason will become clear soon.
{\it Sherman identity}, as we prefer to call it, for short, 
 is a special non trivial case of
another identity 
 called {\it  Feynman identity} 
 first conjectured
by Richard Feynman. 
This identity
relates
the Euler polynomial of a graph to a formal product over the 
classes of equivalence of closed nonperiodic paths with no backtracking in the graph and it is
an important ingredient in a combinatorial formulation of
the Ising model in two dimensions much studied in physics. 
In \cite{sherman} S. Sherman proved 
Feynman identity for planar and toroidal graphs and
 recently this identity was proved in great generality by M. Loebl in \cite{loebl} and D. Cimasoni in \cite{cimasoni}.

Sherman 
compared equation  \eqref{A} with 
the multivariante Witt  identity \cite{witt}:
\begin{equation}  \label{B}
\prod_{m_{1},...,m_{R} \geq 0} (1-z_{1}^{m_{1}}...z_{R}^{m_{R}})^{{\cal M}(m_{1},...,m_{R})}
=1-\sum_{i=1}^{R} z_{i}
\end{equation}
\begin{equation}  \label{C}
{\cal M}(m_{1},...,m_{R}) = \sum_{g|m_{1},...,m_{R}} \frac{\mu(g)}{g} 
\frac{(\frac{N}{g})!}{(\frac{N}{g})(\frac{m_{1}}{g})!...(\frac{m_{R}}{g})!}
\end{equation}
where $N=m_{1}+...+m_{R} > 0$, $\mu$ is the 
M\"obius function defined by the rules: a) $\mu(+1)=+1$, b) $\mu(g)=0$, for
$g=p_{1}^{e_{1}}...p_{q}^{e_{q}}$, $p_{1},...,p_{q}$ primes, and any $e_{i}>1$,
c) $\mu(p_{1}...p_{q})=(-1)^{q}$.
The summation runs over all the common divisors
of $m_{1},...,m_{R}$.

Originally, Witt identity appeared associated with Lie algebras.
In this context the formula gives the dimensions of
 the homogeneous 
subspaces of a finitely generated free Lie algebra $L$. If $L(m_{1}, ...,m_{R})$ is the subspace of $L$ generated by
all homogeneous elements of multidegree $(m_{1}, ...,m_{R})$, then $dim L= {\cal M}$. However,
formula \eqref{C} has many applications in combinatorics as well \cite{moree}. Specially relevant  is that 
${\cal M}$ 
can be interpreted as the number of
equivalence classes of closed non periodic paths which traverse 
counterclockwisely the edges of $G_{R}$, 
the same graph associated to Sherman identity  \eqref{A}. This property is stated in \cite{ssherman} without a proof but
this combinatorial interpretation of Witt formula can be understood
 reinterpreted as a coloring
problem of a necklace with $N$ beads with colors chosen out of a set of $R$ colors such that the coloured beads form a
nonperiodic configuration. In another words, ${\cal M}(m_{1},...,m_{R})$ is the number of nonperiodic coloured necklaces composed of $m_{i}$ ocurrences of  the color $i$, $i=1,...,R$.

In \cite{ssherman} Sherman called attention to this association of both identities \eqref{A}
and \eqref{B}
to paths in the same graph which
motivated him to consider the 
problem of finding a relation of \eqref{A} to Lie algebras.
Interpreting \eqref{A} in algebraic terms means
to relate 
the exponents $N_{\pm}$ to some Lie algebraic data.

An investigation of Sherman's problem was initiated  in \cite{costa} and \cite{costav} and a solution obtained for the univariate case of  identity \eqref{A}.
In the present paper we solve the problem in the multivariate formal  case
which requires important improvements.
The counting method developed in \cite{costa} and \cite{costav} is based on a sign 
formula for a path given in 
terms of data encoded in the 
word representation for the path. 
It played a crucial role in getting formulas for
$N_{\pm}$ in the univariate case. However, 
the counting method based on this sign formula is
complicated. In the present paper we make improvements in 
the counting method in order to apply it
to the multivariate case without depending too much on the sign formula. 
The formula is used here only to prove
a simple Lemma. 

In \cite{kang} S-J. Kang and  M-H. Kim  derived  dimension formulas for the homogeneous spaces of general free graded Lie algebras.
We use some of their results to solve Sherman's problem. At the same time our results give a combinatorial 
realization for some of theirs  in terms of paths in a graph.

The paper is organized as follows. In section \ref{sec : Pre}, we recall
the word representation of 
a path  and some basic definitions. A basic Lemma about the distribution of signs in the set of
words of given length is proved.
In section \ref{sec : Count}, we compute formulas for the numbers of equivalence classes
of closed nonperiodic paths of given length. The first of these generalizes  Witt formula in the sense that it counts paths that traverse the edges  of the graph  in all directions (and no backtrackings). The other formulas give the exponents 
in Sherman's  identity \eqref{A}. We also interpret these formulas in terms of a colouring problem. Sherman's problem, that is, 
to give an algebraic meaning to the exponents in \eqref{A} is solved in section
\ref{sec : alg}.

\section{Preliminaries} \label{sec : Pre}

A path in $G_{R}$ is an ordered sequence of the edges which does not 
necessarily respect their orientation. A path is closed and subjected to the
constraint that it never goes immediately backwards over a  previous edge.

Given $G_{r} \subseteq G_{R}$, denote by $i_{1}, ...,i_{r}$ 
an enumeration of the edges of $G_{r}$ in increasing order.
A closed path of length $N \geq  r$ in $G_{r}$ is best represented by a word
of the form
\begin{equation} \label{D}
D_{j_{1}}^{e_{j_{1}}}D_{j_{2}}^{e_{j_{2}}}
...D_{j_{l}}^{e_{j_{l}}}
\end{equation}
where  $l=r, r+1, ..., N$, $j_{k} \in \{ i_{1}, ...,i_{r} \}$, 
$j_{k} \neq j_{k+1}$, $j_{l} \neq j_{1}$, and 
\begin{equation*}
\sum_{k=1}^{l} \mid e_{j_{k}} \mid=N
\end{equation*}
All edges of $G_{r}$ are traversed by a path so that each 
$i_{k}$ appears at least once in the sequence 
$S_{l}=(j_{1},j_{2},...,j_{l})$.
The order in which the symbols $D_{j}^{e_{j}}$
appear in the word indicates the edges traversed
by $p$ and in which order. 
If the sign of $e_{j}$ is positive (negative)
the path traverses $|e_{j}|$ times edge $j$ following the (opposite of)
edge's orientation.

A word is called periodic if it equals
\begin{equation*}
(D_{j_{1}}^{e_{j_{1}}}D_{j_{2}}^{e_{j_{2}}}
...D_{j_{\alpha}}^{e_{j_{\alpha}}})^{g}
\end{equation*}
for some $g > 1$ and the word between parenthesis is nonperiodic.
Number $g$ is called the period of the word. 
Permuting circularly the symbols $D_{j}^{e_{j}}$ in \eqref{D}
one gets $l$ words that represent the same closed path. 
For example, the word 
$D_{1}^{-2}D_{2}^{+1}D_{1}^{+1}D_{2}^{+3}$ is a circular permutation of
$D_{2}^{+1}D_{1}^{+1}D_{2}^{-3}D_{1}^{-2}$. 
Circular words are taken to be
equivalent
because they represent the
same closed path. Although this is also true for a word
and its inversion
\begin{equation*}
D_{j_{l}}^{-e_{j_{l}}}...D_{j_{1}}^{-e_{j_{1}}}
\end{equation*}
they are not taken equivalent here. This is the reason for the exponent $2$ on the right side of \eqref{A} as originally in \cite{ssherman}.

In section \ref {sec : Count} we consider signed  paths. The sign of a path is given by the formula 
\begin{equation} \label{E}
sign(p)=(-1)^{1+n(p)}
\end{equation}
where 
$n(p)$ is the number of integral revolutions of the tangent vector of $p$.
From this definition it follows that
if $p=(h)^{g}$ is a periodic path with odd
period $g$, then
$sign(p)=sign(h)$. If $g$ is even, $sign(p)=-1$.
The sign of a path  can be computed from its word representation \eqref{D} using 
the formula \cite{costav}
\begin{equation} \label{F}
(-1)^{N+l+T+s+1}
\end{equation}
where $T$ is the number of subsequences in the decomposition of
$S_{l}$ into subsequences (see \cite{costav} for definition and example of a
decomposition)
and $s$ is 
the number of negative exponents in \eqref{D}. It follows from the previous sign formulas
that periodic words with even period have negative sign.
\\
\\
The following Lemma is important in the proof of several results in section \ref {sec : Count}.
It was assumed in \cite{costa} and \cite{costav}
without a proof.

\begin{lem}\label{l2.1} 
Given $G_{r} \subseteq  G_{R}$, consider all paths that traverse each edge of $G_{r}$ at least once (no backtracking allowed) and the set of all representative words 
(periodic or not, circular permutations and inversions included) 
of fixed length $N \geq r >1$. 
Then, half of the words has positive sign and the other half has negative sign.
\end{lem}

\noindent{\bf Proof:} 
It suffices to consider  
the subset of words associated to a fixed sequence 
$S_{l} = (j_{1},j_{2},...,j_{l})$. For this sequence the numbers $N$, $l$ 
and $T$ are fixed. The words with these numbers have signs which
depend only on $s \in \{0,1,2,...,l\}$. For
$N+l+T$ even, the sign of a word is
$(-1)^{s+1}$.
If $l=2k$ there are for each odd value of $s$
\begin{equation*}
\left( \begin{array}{c}
                   2k \\ s
\end{array}\right)
\end{equation*}
words with positive sign. Summing over the odd values of $s$ we get the total 
number of $2^{2k-1}$ words with positive sign. Summing over the even values of 
$s$ we get the same number of words with negative sign. 
If $l=2k+1$ a similar counting gives $2^{2k}$ words
with positive (negative) signs. The case $N+l+T$ odd is analogous. 
$\Box$

\section{Counting paths in $G_{r}$ } \label{sec : Count}

Fix a subgraph $G_{r} \subseteq G_{R}$. 
Call $\theta_{\pm} (m_{i_{1}},...m_{i_{r}})$ the number of equivalence classes of closed nonperiodic paths of length $N \geq r$
and $\pm$ signs
that traverse
$m_{i_{1}}$ times edge $i_{1}$,..., $m_{i_{r}}$ times edge $i_{r}$ of 
$G_{r}$, $m_{i_{j}}>0$, $\forall j=1,...,r$, with no backtracks, $m_{i_{1}}+...+m_{i_{r}}=N$ and zero times the edges in $G_{R} \backslash G_{r}$.  In this section we derive formulas for $\theta
:=\theta_{+} + \theta_{-}$ and $\theta_{\pm}$. Notice that $\theta_{\pm}$ is just another name for the exponents $N_{\pm}$
in (1.1) showing only the nonzero entries in $N_{\pm}$.

Firstly, we compute $\theta$.
In the case $r=1$, a path with
$m_{i}>1$ is periodic. 
The non periodic ones are two, the path
with length $N=1$ and its inversion so that $\theta(m_{i})=0$ if
$m_{i}>1$ and $\theta(m_{i})=2$, if $m_{i}=1$. In the other cases,
$\theta$ is given next.

\begin{thm}\label{th3.1}
For $r=2$, define
\begin{equation} \label{e1}
{\cal F}
\left( \frac{m_{i_{1}}}{g},\frac{m_{i_{2}}}{g}
 \right) =
\sum_{a=1}^{\frac{M}{g}}
\frac{2^{2a}}{a}
\left( \begin{array}{c}
                   \frac{m_{i_{1}}}{g}-1 \\ a-1
\end{array}\right)
\left( \begin{array}{c}
                   \frac{m_{i_{2}}}{g}-1 \\ a-1
\end{array}\right)
\end{equation}
where $M=min\{m_{i_{1}},m_{i_{2}}\}$ and,
if $r \geq  3$,
\begin{equation} \label{e2}
{\cal F}
\left( \frac{m_{i_{1}}}{g},...,
\frac{m_{i_{r}}}{g}
 \right) =
\sum_{a=r}^{\frac{N}{g}}
\frac{2^{a}}{a}
\sum_{ \{ S_{a} \}} \prod_{c=1}^{r}
\left( \begin{array}{c}
                   \frac{m_{i_{c}}}{g}-1 \\ t_{i_{c}}-1
\end{array}\right)
\end{equation}
where
$\{ S_{a} \}$ is the set of sequences $(j_{1},...,j_{a})$ such that
$j_{k} \in \{i_{1},...,i_{r}\}$ and
$j_{k} \neq j_{k+1}$, $j_{a} \neq j_{1}$. Number $t_{i_{c}}$ counts 
how many times edge $i_{c}$ occurs in a sequence $S_{a}$.
Use is made  of the convention that the combination symbol in \eqref{e2}
is zero whenever
$t_{i_{c}} > \frac{m_{i_{c}}}{g}$. Then, 
\begin{equation} \label{e3}
\theta (m_{i_{1}},...,m_{i_{r}})=\sum_{g \mid m_{i_{1}},...m_{i_{r}} }
\frac{\mu(g)}{g} {\cal F} \left( \frac{m_{i_{1}}}{g},...,
\frac{m_{i_{r}}}{g}
 \right)
\end{equation}
The summation is over all the common 
divisors $g$ of $m_{i_{1}},...,m_{i_{r}}$, and
$\mu(g)$ is the M\"obius function.
\end{thm}

\noindent{\bf Proof:}
The number 
${\cal K}(l,m_{i_{1}},...m_{i_{r}})$ of 
words with the same values of $m_{i_{1}},...,m_{i_{r}}$ and
$l \in \{r, r+1,...,N \}$
is given by
\begin{equation*}
{\cal K}(l,m_{i_{1}},...,m_{i_{r}})=2^{l} \sum_{ \{S_{l}\} }
\prod_{c=1}^{r}
\left( \begin{array}{c}
                   m_{i_{c}}-1 \\ n_{i_{c}}-1
\end{array}\right)
\end{equation*}
Let's explain this formula a bit.
Number $n_{i_{c}}$ counts the number of ocurrences of 
edge $i_{c}$ in a sequence $S_{l}=(j_{1},...,j_{l})$.
The combination symbol  counts
the number of unrestricted partitions of $m_{i_{c}}$ into $n_{i_{c}}$ nonzero
positive parts \cite{andrews}
so that the product  times $2^{l}$ 
(there are
$2^{l}$ ways of assigning $+$ and $-$ signs to the exponents in \eqref{D})
gives the total number of words representing 
paths which traverse  $m_{i_{1}}$ times edge $i_{1}$,..., $m_{i_{r}}$ 
times edge $i_{r}$ of 
$G_{r} \subseteq G_{R}$ in all possible ways. 
Then, one sums over all 
sequences $S_{l}$ with
the convention that a combination symbol equals zero whenever 
$m < n$.

In the set of ${\cal K}(l,m_{i_{1}},...,m_{i_{r}})$ words there 
is the subset of nonperiodic words plus 
their
circular permutations and inversions 
and the subset of periodic words if any
whose periods are the common divisors of $l$, and $m_{i_{1}},...,m_{i_{r}}$
plus their 
circular permutations and inversions. 
Denote by  $\overline{{\cal K}(l,m_{i_{1}},...,m_{i_{r}})}$ the number 
of elements in the former set.
The words with period $g$ are of the form
\begin{equation*}
(D_{k_{1}}^{e_{k_{1}}}D_{k_{2}}^{e_{k_{2}}}
...D_{k_{\alpha}}^{e_{k_{\alpha}}})^{g}
\end{equation*}
where $\alpha= l/g$, and $D_{k_{1}}^{e_{k_{1}}}D_{k_{2}}^{e_{k_{2}}}
...D_{k_{\alpha}}^{e_{k_{\alpha}}}$
is nonperiodic  
so that
the number of periodic words with period $g$ 
plus their circular permutations and inversions
is given by $\overline{{\cal K}(l/g,m_{i_{1}}/g,...,m_{i_{r}}/g)}$. Therefore,
\begin{equation*}
{\cal K}(l,m_{i_{1}},...,m_{i_{r}})=\sum_{g \mid l,k,m_{i_{1}},...,m_{i_{r}}} 
\overline{{\cal K} \left( \frac{l}{g},\frac{m_{i_{1}}}{g},..., 
\frac{m_{i_{r}}}{g}\right)}
\end{equation*}
The summation is over all the common divisors  $g$ of $l,m_{i_{1}},...,
m_{i_{r}}$.

Applying M\"obius
inversion formula \cite{apostol} it follows that
\begin{equation} \label{e4}
\overline{{\cal K}(l,m_{i_{1}},...,m_{i_{r}})}=
\sum_{g \mid (l,m_{i_{1}},...,m_{i_{r}})}
\mu(g){\cal K} \left( \frac{l}{g}, \frac{m_{i_{1}}}{g},..., 
\frac{m_{i_{r}}}{g} \right)
\end{equation}
where $\mu$ is the M\"obius function. 
To eliminate
circular permutations divide \eqref{e4} by $l$. 
Summing over all possible values of $l$
one gets a formula for the number 
$\theta(m_{i_{1}},...m_{i_{r}})$:
\begin{equation} \label{e5}
\theta(m_{i_{1}},...,m_{i_{r}})=
\sum_{l=r}^{N} \frac{\overline{{\cal K}(l,m_{i_{1}},...,m_{i_{r}})}}{l}
\end{equation}
Upon substitution of \eqref{e4} into \eqref{e5} one gets, for the case $r\geq 3$,
\begin{equation*}
\theta(m_{i_{1}},...,m_{i_{r}})= \sum_{l=r}^{N} \frac{1}{l}
\sum_{g \mid (l,m_{i_{1}},...,m_{i_{r}})}
\mu(g)
2^{\frac{l}{g}} \sum_{ \{S_{\frac{l}{g}}\} }
\prod_{c=1}^{r}
\left( \begin{array}{c}
                   \frac{m_{i_{c}}}{g}-1 \\ \frac{n_{i_{c}}}{g}-1
\end{array}\right)
\end{equation*}
Proceed now as follows.
For a given common divisor $g$ of $m_{i_{1}},...,m_{i_{r}}$, 
sum over all values of $l$ 
which are multiple of $g$. Then, sum over all possible divisors of 
$m_{i_{1}},...,m_{i_{r}}$.
Write $l=ag$, and 
$n=tg$. In the case $r \geq 3$ one has $r/g \leq a \leq N/g$ but unless 
$g=1$ it is not admissible to have $a<r$ because all $r$ edges of the 
graph should be 
traversed. For this reason, $ r \leq a \leq N/g$.
Result (3.2) follows. In the case $r=2$, 
$l$ is even and, for each $l$, only sequences of the form 
$(i_{1},i_{2},...,i_{1},i_{2})$ with 
$n_{i_{1}}=n_{i_{2}}=\frac{l}{2}$ are possible. Put $l=2a$, $a=1,2,...,
M=min\{m_{1},m_{2}\}$ to get \eqref{e1}.
$\Box$
\\
\\
\noindent{\bf Example 1.} From \eqref{e1}, ${\cal F}(1,1)={\cal
F}(1,2)={\cal F}(2,1)={\cal F}(1,3)={\cal F}(3,1)=4$, ${\cal
F}(2,2)=12$, ${\cal F}(1,4)={\cal F}(4,1)={\cal F}(1,5)={\cal F}(5,1)=4$, ${\cal F}(2,3)={\cal
F}(3,2)=20$, ${\cal F}(2,4)={\cal
 F}(4,2)=28$, ${\cal F}(3,3)=\frac{172}{3}$. From \eqref{e3}, $\theta(1,1)=\theta(1,2)=\theta(2,1)=\theta(1,3)=\theta(3,1)=\theta(1,4)=\theta(4,1)=\theta(1,5)=\theta(5,1)=4$,
$\theta(2,2)=10$, $\theta(2,3)=\theta(3,2)=20$, $\theta(3,3)=56$.
\\
\\
\noindent{\bf Example 2.}
 From \eqref{e2},
${\cal F}(1,1,1)=16$, ${\cal
F}(1,1,2)={\cal F}(1,2,1)={\cal F}(2,1,1)=32$, ${\cal F}(1,2,2)={\cal F}(2,1,2)={\cal F}(2,2,1)=112$, ${\cal F}(1,1,3)=
{\cal F}(1,3,1)={\cal F}(3,1,1)=48$, ${\cal
F}(1,1,4)={\cal F}(1,4,1)={\cal F}(4,1,1)=64$, ${\cal F}(1,2,3)={\cal F}(3,1,2)={\cal F}(2,3,1)={\cal F}(3,2,1)=
{\cal F}(1,3,2)={\cal F}(2,1,3)=256$, ${\cal F}(2,2,2)=1056$. From \eqref{e3},
$\theta(1,1,1)=16$, $\theta(1,1,2)=\theta(2,1,1)=\theta(1,2,1)=32$, $\theta(1,2,2)=
\theta(2,1,2)=\theta(2,2,1)=112$, $\theta(1,1,3)=\theta(3,1,1)=\theta(1,3,1)=48$, $\theta(1,1,4)=\theta(4,1,1)=
\theta(1,4,1)=64$,
$\theta(1,2,3)=\theta(3,1,2)=\theta(2,3,1)=\theta(3,2,1)=\theta(1,3,2)=\theta(2,1,3)=256$,
$\theta(2,2,2)=1048$.
\\
\\
\noindent{\bf Remarks.} 
a) Notice that  $\theta$, likewise Witt formula, is given in terms of M\"obius function. However,  formula \eqref{e3}
counts closed nonperiodic paths traversing the edges of $G_{R}$ in all directions (and no backtracking)
and in that sense generalizes Witt formula. Also, our formula has an algebraic meaning of a dimension. See section \ref{sec : alg}.

\noindent b) If $m_{i_{1}}, ..., m_{i_{r}}$ are coprime, ${\cal F} = \theta$. Otherwise, ${\cal F}$ can be rational. For instance,
${\cal F}(3,3)=172/3$. But ${\cal F}' :=N {\cal F}$, $N=m_{i_{1}}+ ... +m_{i_{r}}$, is always a positive integer which counts the number of all words of length $N$.  For example, in the case $N=4$, $m_{1}=m_{2}=2$, ${\cal F}'=48$. The words are $D_{1}^{\pm 2} D_{2}^{\pm 2}$, $D_{1}^{-1} D_{2}^{+1}D_{1}^{+1} D_{2}^{+1}$, $D_{1}^{+1} D_{2}^{-1}D_{1}^{+1} D_{2}^{+1}$,  $D_{1}^{-1} D_{2}^{-1}D_{1}^{+1} D_{2}^{+1}$,  $D_{1}^{-1} D_{2}^{+1}D_{1}^{+1} D_{2}^{-1}$, 
$D_{1}^{-1} D_{2}^{-1}D_{1}^{-1} D_{2}^{+1}$, and $D_{1}^{-1} D_{2}^{-1}D_{1}^{+1} D_{2}^{-1}$, plus four circular permutations for each of them, and the four periodic words $(D_{1}^{\pm 1} D_{2}^{\pm 1})^{2}$ plus two circular permutations for each.

 In terms of ${\cal F}'$,
\begin{equation*}
\theta (m_{i_{1}},...,m_{i_{r}})=\frac{1}{N} \sum_{g \mid m_{i_{1}},...m_{i_{r}} }
\mu(g) {\cal F}' ( \frac{m_{i_{1}}}{g},...,
\frac{m_{i_{r}}}{g}
 )
\end{equation*}
Although the M\"obius function is negative for some divisors $g$, nevertheless the right hand side is always a positive number because ${\cal F}' ( \frac{m_{i_{1}}}{g},...,
\frac{m_{i_{r}}}{g}) $ counts words in a subset of the words counted by ${\cal F}' ( m_{i_{1}},...,
m_{i_{r}})$.

c) Given a circular necklace with $N$ beads
consider the problem of counting inequivalent nonperiodic colourings of these beads with $2r$ colors $\{ c_{i}, \overline{ c}_{i} \}$, $i=1,...,r$,  with $m_{i}$ occurrences of the index $i$, $N=\sum m_i$,  with the restriction that no two colors $c_i$ and $\overline{c_i}$  (same index) occur  adjacent  in a colouring.  Now, consider an oriented graph with $r$ loops hooked to a single vertex. Each loop edge corresponds to a color $c_i$. A nonperiodic closed nobacktracking path of length $N$ in the graph  corresponds to a colouring and a color $\overline c_i$ corresponds to  an edge being traversed in the opposite orientation. The presence of a single vertex in the graph reflects the fact that adjacent to a bead with, say,  color $c_{i}$  any other with distinct index may follow. The number of inequivalent colourings is given by $\theta$.
\\
\\
As a basic test of our counting ideas, we prove Sherman's statement in \cite{ssherman} relating Witt formula to paths in $G_{R}$:

\begin{prop}\label{pro3.2} 
Relative to graph $G_{R}$,
formula (1.2) gives the number ${\cal M}$ of equivalence classes
of closed non periodic paths of length $N>0$ which traverse
counterclockwisely $m_i \geq 0$ times edge $i$, $i=1,2,...,R$,
$m_{1}+...+m_{R}=N$.
\end{prop}

\noindent{\bf Proof:} Denote by $m_{i_{1}}$,...,$m_{i_{r}}$, $r
\leq R$, the non zero entries in ${\cal M} (m_{1},...,m_{R})$ which we call ${\cal M}_{r} (m_{i_{1}},...,m_{i_{r}})$. Words representing
counterclockwise paths have positive exponents so that the factors
$2^{2a}$ and $2^{a}$ in formulas \eqref{e1} and \eqref{e2} are not
needed, hence,
\begin{equation} \label{e6}
{\cal M}_{r} (m_{i_{1}},...,m_{i_{r}}) = \sum_{g \mid
m_{i_{1}},...,m_{i_{r}} } \frac{\mu(g)}{g} {\cal F}_{c} \left(
\frac{m_{i_{1}}}{g},..., \frac{m_{i_{r}}}{g}
 \right)
\end{equation}
where
\begin{equation} \label{e7}
{\cal F}_{c} \left( \frac{m_{i_{1}}}{g},\frac{m_{i_{2}}}{g}
 \right) =
\sum_{a=1}^{\frac{M}{g}} \frac{1}{a} \left( \begin{array}{c}
                   \frac{m_{i_{1}}}{g}-1 \\ a-1
\end{array}\right)
\left( \begin{array}{c}
                   \frac{m_{i_{2}}}{g}-1 \\ a-1
\end{array}\right)
\end{equation}
with $M=min\{ m_{i_{1}},m_{i_{2}} \}$, if $r=2$; and
\begin{equation} \label{e8}
{\cal F}_{c} \left( \frac{m_{i_{1}}}{g},... \frac{m_{i_{r}}}{g}
 \right) =
\sum_{a=r}^{\frac{N}{g}} \frac{1}{a}
\sum_{ \{ S_{a} \}} \prod_{c=1}^{r}
\left( \begin{array}{c}
                   \frac{m_{i_{c}}}{g}-1 \\ t_{i_{c}}-1
\end{array}\right)
\end{equation}
if $r \geq 3$. In the case $r=2$ suppose $m_{i_{1}} \leq
m_{i_{2}}$. Using formula \eqref{5.3} (with $l=2$), section \ref{sec : Lem}, it follows that 
\begin{eqnarray*}
\sum_{a=1}^{\frac{m_{i_{1}}}{g}} \frac{1}{a} \left(
\begin{array}{c}
                   \frac{m_{i_{1}}}{g}-1 \\ a-1
\end{array}\right)
\left( \begin{array}{c}
                   \frac{m_{i_{2}}}{g}-1 \\ a-1
\end{array}\right)
& = & \frac{g}{m_{i_{2}}} \left( \begin{array}{c}
                   \frac{m_{i_{1}}}{g}+\frac{m_{i_{2}}}{g}-1 \\
\frac{m_{i_{1}}}{g}
\end{array}\right) \\
& = & \frac{(\frac{N}{g})!} {(\frac{N}{g})(\frac{m_{i_{1}}}{g})!
(\frac{m_{i_{2}}}{g} )!}
\end{eqnarray*}
Similarly, if $m_{i_{2}} \leq m_{i_{1}}$. In the case $r \geq 3$
define $I$,
\begin{equation} \label{e9}
I= \sum_{ \begin{array}{c}
      \footnotesize{m_{i} >0}\\
       \footnotesize{m_{i_{1}}+...+m_{i_{r}}=N}
       \end{array}}
       {\cal F}_{c} \left( \frac{m_{i_{1}}}{g},...,\frac{m_{i_{r}}}{g} \right)
\end{equation}
Upon substitution of \eqref{e8} into \eqref{e9} and exchanging the
summation symbols, we get
\begin{equation*}
I=\sum_{a=r}^{\frac{N}{g}} \frac{1}{a} \sum_{\{S_{a}\}} \sum_{
\begin{array}{c}
      \footnotesize{m_{i} >0}\\
       \footnotesize{m_{i_{1}}+...+m_{i_{r}}=N}
       \end{array}}
        \prod_{c=1}^{r}
\left( \begin{array}{c}
                   \frac{m_{i_{c}}}{g}-1 \\ t_{i_{c}}-1
\end{array}\right)
\end{equation*}

\noindent Applying Lemma \ref{l5.2}, section \ref{sec : Lem},
\begin{equation*}
I= \sum_{a=r}^{\frac{N}{g}} \frac{1}{a} \sum_{\{S_{a}\}}
\left( \begin{array}{c}
                   \frac{N}{g}-1 \\ a-1
\end{array}\right)
=
\sum_{a=r}^{\frac{N}{g}} \frac{1}{a}
\left( \begin{array}{c}
                   \frac{N}{g}-1 \\ a-1
\end{array}\right)
rw_{r}(a)
\end{equation*}
where
\begin{equation*}
rw(a)= \sum_{j=1}^{r} (-1)^{r+j}
\left( \begin{array}{c}
                   r \\ j
\end{array}\right) (j-1)^{a}+(-1)^{a+r}
\end{equation*}
is the number of sequences in $\{S_{a}\}$ \cite{costav}.  Using that
\begin{equation*}
\sum_{a=r}^{\frac{N}{g}} \frac{1}{a}
\left( \begin{array}{c}
                   \frac{N}{g}-1 \\ a-1
\end{array}\right) (j-1)^{a} =\frac{g}{N} (j^{\frac{N}{g}}-1)
\end{equation*}
and
\begin{equation*}
\sum_{a=r}^{\frac{N}{g}} \frac{1}{a}
\left( \begin{array}{c}
                   \frac{N}{g}-1 \\ a-1
\end{array}\right) (-1)^{a+r} = (-1)^{r+1} \frac{g}{N}
\end{equation*}
we get
\begin{equation} \label{e10}
I = \frac{g}{N} \sum_{j=1}^{r} (-1)^{r+j}
\left( \begin{array}{c}
                   r \\ j
\end{array}\right)
j^{\frac{N}{g}}
\end{equation}
Stirling numbers $S(\frac{N}{g}, r)$ of
second kind are given by the formula \cite{chuan}
\begin{equation} \label{e11}
S\left( \frac{N}{g}, r \right)=\frac{1}{r!} \sum_{k=0}^{r}(-1)^{k}
\left( \begin{array}{c}
                   r \\ k
\end{array}\right) (r-k)^{\frac{N}{g}} =
\frac{1}{r!} \sum_{j=0}^{r}(-1)^{r+j}
\left( \begin{array}{c}
                   r \\ j
\end{array}\right) j^{\frac{N}{g}}
\end{equation}
so that
\begin{equation} \label{e12}
I= r! \frac{g}{N} S \left( \frac{N}{g},r \right) 
\end{equation}
Stirling numbers have the property that
\begin{equation} \label{e13}
\sum_{ \begin{array}{c}
      \footnotesize{m_{i} >0}\\
       \footnotesize{m_{i_{1}}+...+m_{i_{r}}=N}
       \end{array}} \frac{(\frac{N}{g})!}
{ (\frac{m_{i_{1}}}{g})!...(\frac{m_{i_{r}}}{g})! }= r!
S \left( \frac{N}{g},r \right) 
\end{equation}
Comparing relations \eqref{e12}, \eqref{e13} and \eqref{e9},
\begin{equation} \label{e14}
{\cal F}_{c} \left( \frac{m_{i_{1}}}{g},...,\frac{m_{i_{r}}}{g} \right) =
\frac{g}{N} \frac{(\frac{N}{g})!}{(\frac{m_{1}}{g})!...
(\frac{m_{r}}{g})!} 
\end{equation}
Upon substitution of \eqref{e14} into \eqref{e6} the result follows.
\\
\\
In the sequel we compute formulas for $\theta_{+}$ and $\theta_{-}$.

\begin{thm} \label{th3.3}
Suppose any of the following conditions is 
satisfied: 
\noindent{\bf (a)} $N=m_{i_{1}}+...+m_{i_{r}}< 2r$;
\noindent{\bf (b)} $m_{i_{1}}, ...,m_{i_{r}}$ are coprime;
\noindent{\bf (c)} $m_{i_{1}}, ...,m_{i_{r}}$ are not all odd nor even;
\noindent{\bf (d)} $m_{i_{1}}, ...,m_{i_{r}}$ are all odd.
\noindent Then,
\begin{equation} \label{e15}
\theta_{-}(m_{i_{1}}, ..., m_{i_{r}})=
\theta_{+}(m_{i_{1}},...,m_{i_{r}})
\end{equation} 
\end{thm}
\noindent{\bf Proof:} Similar to Theorem 1 in \cite{costa} using Lemma 1. $\Box$.

The case where $m_{i_{1}},...,m_{i_{r}}$ are all even numbers is given in the 
next theorem.

\begin{thm} \label{th3.4}
The number $\theta_{+} (m_{i_{1}},...m_{i_{r}})$ 
is given by
\begin{equation} \label{e16}
\theta_{+} (m_{i_{1}},...,m_{i_{r}})=\sum_{odd \hspace{1mm} g \mid 
m_{i_{1}},...,m_{i_{r}}}
\frac{\mu(g)}{g} {\cal G} \left( \frac{m_{i_{1}}}{g},...,
\frac{m_{i_{r}}}{g}
 \right)
\end{equation}
where
the summation is over all the common 
odd divisors of $m_{i_{1}},...,m_{i_{r}}$, and 
${\cal G} = \frac{\cal F}{2}$
with $\cal F$ 
as in \eqref{e1} and \eqref{e2}. Suppose $m_{i_{1}},...,m_{i_{r}}$ are all even numbers.
Then,
\begin{equation} \label{e17}
\theta_{-}(m_{i_{1}}, ..., m_{i_{r}})=\theta_{+}
(m_{i_{1}},...,m_{i_{r}}) - \theta_{+} \left(
\frac{m_{i_{1}}}{2},...,\frac{m_{i_{r}}}{2}, \right)
\end{equation}
\end{thm}
\noindent{\bf Proof:}
First, suppose that all common divisors of $m_{i_{1}},...,m_{i_{r}}$
are odd numbers. In this case,
\begin{equation*}
\theta (m_{i_{1}},...,m_{i_{r}} ) 
= \sum_{odd \hspace{1mm} g \mid m_{i_{1}},...,m_{i_{r}}}
\frac{\mu(g)}{g} {\cal F} \left(  \frac{m_{i_{1}}}{g},...,\frac{m_{i_{r}}}{g}
 \right)
\end{equation*}
Since $\theta = \theta_{+} + \theta_{-}$ and $\theta_{+}=\theta_{-}$ 
(Theorem \ref{th3.2})
it follows that $\theta = 2 \theta_{+}$, hence,
\begin{equation} \label{e18}
\theta_{+} = \frac{1}{2}
\sum_{odd \hspace{1mm} g \mid m_{i_{1}},...,m_{i_{r}}  }
\frac{\mu(g)}{g} {\cal F} \left( \frac{m_{i_{1}}}{g},...,\frac{m_{i_{r}}}{g}   \right)
\end{equation}
If the numbers $m_{i_{1}},...,m_{i_{r}}$  are all even then again $\theta_{+}$
is given  by \eqref{e18} for in this case the $m_{i}$'s have common 
divisors which are even numbers but since periodic words with even 
period have negative sign, hence, only the odd divisors are relevant to 
get
$\theta_{+}$. The reason why one should have the factor $1/2$ 
is that by Lemma 1 when one considers the set of all possible words 
representing paths of
a given length which traverse $m_{i_{1}},...,m_{i_{r}}$ times
the edges of $G_{r}$, half of them have
positive sign and the other half have negative sign. To account for the 
positive half one needs the factor $1/2$.
Let's now compute $\theta_{-}$ in the even case. Write
\begin{eqnarray*}
\theta & = & \sum_{odd \hspace{1mm} g \mid m_{i_{1}},...,m_{i_{r}}}
\frac{\mu(g)}{g} {\cal F}  +
\sum_{even \hspace{1mm} g \mid m_{i_{1}},...,m_{i_{r}}}
\frac{\mu(g)}{g} {\cal F} \\
& = & \frac{1}{2} \sum_{odd \hspace{1mm} g \mid m_{i_{1}},...,m_{i_{r}}}
\frac{\mu(g)}{g} {\cal F}  +
 \frac{1}{2} \sum_{odd \hspace{1mm} g \mid m_{i_{1}},...,m_{i_{r}}}
\frac{\mu(g)}{g} {\cal F} \\
&  & + \sum_{even \hspace{1mm} g \mid m_{i_{1}},...,m_{i_{r}}}
\frac{\mu(g)}{g} {\cal F} \\
& = & 2 \theta_{+} +
\sum_{even \hspace{1mm} g \mid m_{i_{1}},...,m_{i_{r}}}
\frac{\mu(g)}{g} {\cal F}
\end{eqnarray*}
Using that $\theta=\theta_{+}+\theta_{-}$, it follows that
\begin{equation*}
\theta_{-}= \theta_{+} + \sum_{even \hspace{1mm} g \mid m_{i_{1}},...,m_{i_{r}}}
\frac{\mu(g)}{g} {\cal F}
\end{equation*}
Now, the relevant even divisors are $\{ 2n \}$ where $n$ are the odd common
divisors of  $\{m_{i}\}$. For the other possible divisors if any use that 
$\mu (2^{j} n) = 0$, $j \geq 2$. Using that $\mu (2n)= - \mu (n)$
the summation over the even divisors is equal to 
\begin{equation*}
-\theta_{+} ( \frac{m_{i_{1}}}{2},...,\frac{m_{i_{r}}}{2})
\end{equation*}
proving the result.
$\Box$
\\
\\
\noindent{\bf Remark.}  Likewise $\theta$, the numbers $\theta_{\pm}$ can be interpreted as the number of inequivalent nonperiodic colourings of a circular necklace with $N$ beads.
However, now these colourings are classified as positive or negative according to formula \eqref{F}. It is positive (negative) if the number $N+l+T+s$ is odd (even). In this case, $s$ is the number of $\overline{c}$ colors present in a colouring. Interpret $T$ in terms of the color indices.
\\
\\
{\it \noindent{\bf Definition.} Let $s_{1},...,s_{r}$ be
arbitrary positive integers. Let the number ${\cal P}$ be defined as follows. 
If
$s_{1},...,s_{r}$ are all even numbers,
\begin{equation} \label{e19}
{\cal P}(s_{1},...,s_{r})=\sum_{even \hspace{1mm} g \mid s_{1},...,s_{r}}
\frac{\mu(g)}{g} {\cal G} 
\left( \frac{s_{1}}{g},...,\frac{s_{r}}{g} 
\right)
\end{equation}
Otherwise,
${\cal P}(s_{1},...,s_{r})=0$. Also, define
\begin{equation} \label{e20}
{\cal H}
=
\left\{
\begin{array}{ll}
{\cal G}(s_{1},...,s_{r}) & \mbox{if $s_{1},...,s_{r}$ not all even}\\
 {\cal G}(s_{1},...,s_{r}) - \sum_{k \mid s_{1},...,s_{r}}
\frac{1}{k} {\cal P}(\frac{s_{1}}{k},...,\frac{s_{r}}{k}) &
\mbox{otherwise}
\end{array}
\right.
\end{equation}}

\begin{lem}\label{lem3.5}
\begin{equation} \label{e21}
{\cal P}=\sum_{g \mid s_{1},...,s_{r} } 
\frac{\mu(g)}{g} \left( {\cal G}- {\cal H} \right)
\end{equation}
\end{lem}
\noindent{\bf Proof:} From the above definition, ${\cal G} = {\cal H}$ if
$s_{1},...,s_{r}$ not all even. Otherwise,
\begin{equation*}
{\cal G}-{\cal H}
= \sum_{g \mid s_{1},...,s_{r}}
\frac{1}{g} {\cal P} \left(\frac{s_{1}}{g},...,\frac{s_{r}}{g} \right)
\end{equation*}
Now, apply Lemma 5.1, section \ref{sec : Lem},  to get the result. $\Box$

\begin{thm}\label{thm3.6}
\begin{equation} \label{e22}
\theta_{+}(m_{i_{1}},...,m_{i_{r}})
=
\sum_{g | m_{i_{1}},...,m_{i_{r}}} \frac{\mu(g)}{g} 
{\cal H} \left(  \frac{m_{i_{1}}}{g},...,\frac{m_{i_{r}}}{g}   \right)
\end{equation}
\end{thm}
\noindent{\bf Proof:}
When $m_{i_{1}},...,m_{i_{r}}$ are not all even, their odd divisors 
are the only possible common divisors. In this case, ${\cal P}=0$ and
$$
\theta_{+} = \sum_{odd \hspace{1mm} g \mid m_{i_{1}},...,m_{i_{r}}}
\frac{\mu(g)}{g} {\cal H} 
$$
with ${\cal H} ={\cal G}$.
In the case $m_{i_{1}},...,m_{i_{r}}$ are all even 
the sum over odd divisors of $m_{i_{1}},...,m_{i_{r}}$ can be expressed as
\begin{eqnarray*}
\theta_{+}
& = &
\sum_{odd \hspace{1mm} g | m_{i_{1}},...,m_{i_{r}}} \frac{\mu(g)}{g} 
{\cal G}  \\
& = &
\sum_{ g | m_{i_{1}},...,m_{i_{r}}} \frac{\mu(g)}{g} 
{\cal G}  
- \sum_{even \hspace{1mm} g | m_{i_{1}},...,m_{i_{r}}} \frac{\mu(g)}{g} 
{\cal G} 
\\
& = & \sum_{ g | m_{i_{1}},...,m_{i_{r}}} \frac{\mu(g)}{g} 
{\cal G} 
- {\cal P} \\
& = &
\sum_{ g | m_{i_{1}},...,m_{i_{r}}} \frac{\mu(g)}{g} 
{\cal G} 
- \sum_{ g | m_{i_{1}},...,m_{i_{r}}} \frac{\mu(g)}{g} ( {\cal G}-{\cal H} ) \\
& = & \sum_{ g | m_{i_{1}},...,m_{i_{r}}} \frac{\mu(g)}{g} {\cal H}
\end{eqnarray*}
$\Box$
\\
\\
\noindent{\bf Example 3.}
$\theta_{\pm}(1,1)=\theta_{\pm}(1,2)=\theta_{\pm}(2,1)=\theta_{\pm}(1,3)=\theta_{\pm}(3,1)=\theta_{\pm}(1,4)=
\theta_{\pm}(4,1)=\theta_{\pm}(1,5)=\theta_{\pm}(5,1)=2$,
$\theta_{+}(2,2)=6$, $\theta_{-}(2,2)=4$,
$\theta_{\pm}(2,3)=\theta_{\pm}(3,2)=10$,
$\theta_{+}(2,4)=14$, $\theta_{-}(2,4)=12$, $\theta_{+}(4,2)=14$, $\theta_{-}(4,2)=12$, $\theta_{\pm}(3,3)=28$. 
\\
\\
\noindent{\bf Example 4.}
$\theta_{\pm}(1,1,1)=8$, $\theta_{\pm}(1,1,2)=\theta_{\pm}(2,1,1)=\theta_{\pm}(1,2,1)=16$, $\theta_{\pm}(1,2,2)=
\theta_{\pm}(2,1,2)=\theta_{\pm}(2,2,1)=56$, $\theta_{\pm}(1,1,3)=\theta_{\pm}(3,1,1)=\theta_{\pm}(1,3,1)=24$,$\theta_{\pm}(1,1,4)=\theta_{\pm}(4,1,1)=\theta_{\pm}(1,4,1)=32$,
$\theta_{\pm}(1,2,3)=\theta_{\pm}(3,1,2)=\theta_{\pm}(2,3,1)=\theta_{\pm}(3,2,1)=\theta_{\pm}(1,3,2)=\theta_{\pm}(2,1,3)=128$,
$\theta_{+}(2,2,2)=524$, $\theta_{-}(2,2,2)=516$.

\section{Sherman identity and Lie algebras} \label{sec : alg}

In this section we relate our previous results with Lie algebras
and solve Sherman's problem. The solution is provided by the following proposition by S. -J. Kang and M. -H. Kim  in \cite{kang}.

\begin{prop}\label{pro3.7}
Let 
$V = \bigoplus_{(k_{1}, ..., k_{r}) \in {\bf Z}_{>0}^{r}} 
V_{(k_{1}, ..., k_{r}) }$ 
be a ${\bf Z}_{>0}^{r}$-graded vector
space over ${\bf C}$ with $dim V_{(k_{1}, ...,k_{r}) }= 
d(k_{1}, ..., k_{r}) < \infty$, for all
$(k_{1}, ..., k_{r}) \in {\bf Z}_{>0}^{r}$,
and let 
$L= \bigoplus_{(k_{1}, ..., k_{r}) \in {\bf Z}_{>0}^{r}  } 
L_{(k_{1}, ..., k_{r})}$ be the free Lie algebra
generated by $V$.
Then, the dimensions of the subspaces $L_{(k_{1}, ..., k_{r}) }$
are given by
\begin{equation} \label{4.1}
dim L_{(k_{1}, ..., k_{r})}
= \sum_{g | (k_{1}, ..., k_{r})} \frac{\mu(g)}{g} {\cal W} \left(\frac{k_{1}}{g},
...,\frac{k_{r}}{g}\right)
\end{equation}
where summation is over all common divisors $g$ of $k_{1}, ..., k_{r}$ 
and ${\cal W}$ is given by
\begin{equation} \label{4.2}
{\cal W}(k_{1}, ..., k_{r}) = \sum_{s \in T(k_{1}, ..., k_{r})}
\frac{(|s|-1)!}{s!} \prod_{i_{1},...,i_{r}=1}^{\infty} 
d(i_{1}, ..., i_{r})^{s_{i_{1}, ..., i_{r}}}
\end{equation}
The exponents $s_{i_{1}, ..., i_{r}}$ are the components of $s \in T$,
\begin{multline} \label{4.3}
T(k_{1}, ..., k_{r})=\{s=(s_{i_{1}, ..., i_{r}})|s_{i_{1}, ..., i_{r}}
\in {\cal Z}_{\geq 0}, \\
\sum_{i_{1},...,i_{r}=1}^{\infty}
s_{i_{1}, ..., i_{r}} (i_{1},...,i_{r})=
(k_{1}, ..., k_{r})\},
\end{multline}
and
\begin{equation} \label{4.4}
|s|=\sum_{i_{1},...,i_{r}=1}^{\infty}
s_{i_{1}, ..., i_{r}}, \hspace{5mm} s!= \prod_{i_{1},...,i_{r}=1}^{\infty}
s_{i_{1}, ..., i_{r}}!
\end{equation}
Moreover, the numbers $dim L_{(k_{1}, ..., k_{r}) } $ satisfy
\begin{equation} \label{4.5}
\prod_{k_{1},...,k_{r}=1}^{\infty} 
(1-z_{1}^{k_{1}}...z_{r}^{k_{r}})^{dim L_{(k_{1}, ..., k_{r})}}= 
1 - f(z_{1},...,z_{r})
\end{equation}
where
\begin{equation} \label{4.6}
f(z_{1},...,z_{r}) := \sum_{k_{1}, ..., k_{r}=1}^{\infty} d(k_{1}, ..., k_{r})
z_{1}^{k_{1}}...z_{r}^{k_{r}}
\end{equation}
This function is associated with the generating function of the ${\cal W}$'s,
\begin{equation} \label{4.7}
g(z_{1},...,z_{r}) :=\sum_{k_{1}, ..., k_{r}=1}^{\infty} {\cal W}(k_{1}, ..., k_{r} )
z_{1}^{k_{1}}...z_{r}^{k_{r}}
\end{equation}
by the relation
\begin{equation} \label{4.8}
e^{-g}=  1-f
\end{equation}
\end{prop}
$\Box$
\\
\\
Identity \eqref{4.5}  is a consequence of the famous {\it Poincar\'e-Birkhoff-Witt theorem} for the free Lie algebra.  Computation of the formal logarithm of the left hand side  of \eqref{4.5} and its expansion gives that the infinite product equals the exponential in \eqref{4.8}. Raise both members of \eqref{4.5} to the power $-1$, compute  the formal logarithm of both members and expand them. Identification of the coefficients of the same order,  definition \eqref{4.2} and application of  M\"obius inversion gives \eqref{4.1}. See \cite{kang} for  details. In \cite{kang}, \eqref{4.1} is called the {\it generalized Witt formula}, ${\cal W}$ is called  the {\it Witt partition function} and \eqref{4.5} the {\it
generalized Witt 
identity}.

Formulas \eqref{e3} and \eqref{e22} have exactly the form of  \eqref{4.1} with corresponding  {\it Witt partition functions} given by 
${\cal F}$, ${\cal H}$, respectively, so  we will
interpret  $\theta$ and $\theta_{+}$ as giving the dimensions of the 
homogeneous spaces of  graded Lie algebras. 
In each case, the algebra  is generated by 
a graded vector space whose dimensions 
can be computed recursively from \eqref{4.2} as a function of the Witt partition function. However, a general formula can be obtained
from \eqref{4.8} using \eqref{4.6} as the formal Taylor
expansion of $1-e^{-g}$. This gives
\begin{equation} \label{4.9}
d(k_{1},...,k_{r}) =
\frac{1}{k_{1}!...k_{r}!} 
\frac{\partial^{|k|}}{\partial z_{1}^{k_{1}}...
\partial z_{r}^{k_{r}}} (1-e^{-g})|_{z_{1}=...=z_{r}=0}
\end{equation}
with
\begin{equation} \label{4.10}
g(z_{1},...,z_{r}) :=\sum_{k_{1}, ..., k_{r}=1}^{\infty} {\cal W}(k_{1}, ..., k_{r} )
z^{k_{1}}...z^{k_{r}}
\end{equation}
and ${\cal W}= {\cal F}, {\cal H}$  given by \eqref{e1}, \eqref{e2}, \eqref{e20}. Furthemore,
$dim L_{(k_{1}, ..., k_{r})}= \theta, \theta_{+}$ 
given by \eqref{e3}, \eqref{e22}
satisfy the
{\it generalized Witt identity}  \eqref{4.5} with the corresponding dimensions given by \eqref{4.9}.
In fact, an explicit formula for \eqref{4.9} can be derived as follows:

\begin{thm} \label{thm4.2}
A formula for
the numbers $d(k_{1},...,k_{r})$ 
is
\begin{equation} \label{4.11}
d(k_{1},...,k_{r})=\sum_{\lambda =1}^{|k|} (-1)^{\lambda+1}
\sum_{p(\lambda,k)} \prod_{i=1}^{q} 
\frac{[{\cal W}({ l_{i1}}, ..., { l_{ir}} )]^{a_{i}}}{a_{i}!}
\end{equation}
where $|k|=k_{1}+...+k_{r}$, $q=-1+\prod_{i=1}^{r}(k_{i}+1)$,
$p_{\lambda,k}$ is the set of all $a_{i}\in\{0,1,2,...\}$ such that
$\sum_{i=1}^{q}a_{i}=\lambda$, $\sum_{i=1}^{q}a_{i}l_{ij}=k_{j}$, and 
the vectors $l_{i}=(l_{i1},...,l_{ir})$, $l_{ij}$ satisfying 
$0 \leq l_{ij} \leq k_{j}$, $\forall j=1,...,r$, $\forall i=1,...,q$ and 
$\sum_{j=1}^{r}l_{ij}>0$. Set ${\cal W}(l_{i})=0$ if $l_{ij}=0$ for some $j$; otherwise, 
${\cal W}$ is the Witt partition function.
\end{thm}

\noindent{\bf Proof:}
A generalization of
Fa\`a di Bruno's relation due to Constantine 
and Savits in \cite{savits} and \cite{ssavits} gives 
a formula for the $|k|$-th derivative
of the exponential of a function $g(z_{1},...,z_{r})$. 
From this formula and \eqref{4.9}, \eqref{4.11} follows. 
$\Box$.
\\
\\
\noindent{\bf Example 5.}
We compute $d(2,2)$, explicitly. In this case, $k_{1}=
k_{2}=2$, $|k|=4$, $q=8$. The possible vectors $l \leq (2,2)$ are 
$l_{1}=(0,1)$, $l_{2}=(1,0)$, $l_{3} = (1,1)$, $l_{4}=(0,2)$, $l_{5}=(2,0)$,
$l_{6}=(2,1)$, $l_{7}=(1,2)$ and $l_{8}=(2,2)$. Next we give the values of
$a_{1},...,a_{8} \geq 0$ satisfying
\begin{equation*}
\sum_{i=1}^{8} a_{i} = \lambda, \hspace{1cm} 
\sum_{i=1}^{8} a_{i}l_{i}= (2,2)
\end{equation*}
Define the vector $a=(a_{1},...,a_{8})$. The possible $a$'s for each $\lambda$
are as follows. 
For $\lambda =1$, $a=(0,...,0,1)$. For $\lambda =2$, $(0,1,0,0,0,0,1,0)$,
$(0,0,2,0,0,0,0,0)$, $0,0,0,1,1,0,0,0)$, $(1,0,0,0,0,1,0,0)$. For $\lambda=3$,
$(0,2,0,1,0,0,0,0)$, $(2,0,0,0,1,0,0,0)$, $(1,1,1,0,0,0,0,0)$. For $\lambda=4$,
$(2,2,0,0,0,0,0,0)$. We get
\begin{equation*}
d(2,2) = {\cal W}(2,2) - \frac{1}{2} {\cal W}(1,1)^{2}
\end{equation*}
The dimensions up to $d(3,3)$ are:
\begin{eqnarray*}
N=2  \hspace{1cm} d(1,1) &=& {\cal W}(1,1)\\
N=3  \hspace{1cm} d(1,2) &=& {\cal W}(1,2), \hspace{5mm} d(2,1) = {\cal W}(2,1)\\
N=4  \hspace{1cm} d(1,3) &=& {\cal W}(1,3), \hspace{5mm} d(3,1) = {\cal W}(3,1)\\ 
d(2,2) &=& {\cal W}(2,2) - \frac{1}{2} {\cal W}(1,1)^{2}\\
N=5   \hspace{1cm} d(1,4) &=& {\cal W}(1,4), \hspace{5mm} d(4,1) = {\cal W}(4,1)\\
d(2,3) &=& {\cal W}(2,3) - {\cal W}(1,1) {\cal W}(1,2)\\
d(3,2) &=& {\cal W}(3,2) - {\cal W}(1,1) {\cal W}(2,1)\\
N=6 \hspace{1cm} d(1,5) &=& {\cal W}(1,5), \hspace{5mm} d(5,1) = {\cal W}(5,1)\\
d(2,4) &=& {\cal W}(2,4) - {\cal W}(1,1) {\cal W}(1,3)-\frac{1}{2} {\cal W}(1,2)^{2}\\
d(4,2) &=& {\cal W}(4,2) - {\cal W}(1,1) {\cal W}(3,1)-\frac{1}{2} {\cal W}(2,1)^{2}\\
d(3,3) &=& {\cal W}(3,3) - {\cal W}(1,1){\cal W}(2,2)-{\cal W}(1,2){\cal W}(2,1)+\frac{1}{6}{\cal W}(1,1)^{3}
\end{eqnarray*}

\noindent For $r=3$, the dimensions up to $d(2,2,2)$ are:
\begin{eqnarray*}
N=3  \hspace{1cm} d(1,1,1) &=& {\cal W}(1,1,1)\\
N=4  \hspace{1cm} d(1,1,2) &=& {\cal W}(1,1,2), \hspace{2mm} d(1,2,1) = {\cal W}(1,2,1), \hspace{2mm} d(2,1,1)={\cal W}(2,1,1)\\
N=5   \hspace{1cm} d(1,2,2) &=& {\cal W}(1,2,2), \hspace{2mm} d(2,1,2) = {\cal W}(2,1,2), \hspace{2mm} d(2,2,1)={\cal W}(2,2,1)\\
d(1,1,3) &=& {\cal W}(1,1,3), \hspace{2mm} d(1,3,1)={\cal W}(1,3,1), \hspace{2mm} d(3,1,1)={\cal W}(3,1,1)\\
N=6 \hspace{1cm} d(1,1,4) &=& {\cal W}(1,1,4), \hspace{2mm} d(1,4,1) = {\cal W}(1,4,1), \hspace{2mm} d(4,1,1)={\cal W}(4,1,1)\\
d(1,2,3) &=& {\cal W}(1,2,3),\hspace{2mm} d(3,1,2) ={\cal W}(3,1,2), \hspace{2mm} d(2,3,1) ={\cal W}(2,3,1)\\
d(3,2,1) &=& {\cal W}(3,2,1),\hspace{2mm} d(1,3,2) = {\cal W}(1,3,2), \hspace{2mm} d(2,1,3)={\cal W}(2,1,3)\\
d(2,2,2) &=&{\cal W}(2,2,2)-\frac{1}{2} {\cal W}^{2}(1,1,1)\\
\end{eqnarray*}
\noindent{\bf Example 6.} Relative to $\theta$  with ${\cal W}={\cal G}$ and applying data from previous examples  for the case $r=2$ we find
 the dimensions
$d(1,1)=d(1,2)=d(2,1)=d(1,3)=d(3,1)=d(1,4)=d(4,1)=d(2,3)=d(3,2)=d(1,5)
=d(5,1)= d(2,2)= d(2,4)=d(4,2)=d(3,3)=4$.  In the case $r=3$, the dimensions are
$d(1,1,1)=8$, $d(1,1,2)=d(2,1,1)=d(1,2,1)=16$, $d(1,2,2)=d(2,1,2)=d(2,2,1)=56$, $d(1,1,3)=d(3,1,1)=d(1,3,1)=24$,
$d(1,1,4)=d(4,1,1)=d(1,4,1)=32$, $d(1,2,3)=d(3,1,2)=d(2,3,1)=d(3,2,1)=d(1,3,2)=d(2,1,3)=128$, $d(2,2,2)=496$. 
\\
\\
\noindent{\bf Example 7.} Relative to $\theta_{+}$  with ${\cal W}={\cal H}$  we find the dimensions
$d(1,1)=d(1,2)=d(2,1)=d(1,3)=d(3,1)=d(1,4)=d(4,1)=d(2,3)=d(3,2)=d(1,5)
=d(5,1)=2, d(2,2)=5, d(2,4)=d(4,2)=9, d(3,3)=28$ for the case $r=2$. Also,
$d(1,1,1)=8$, $d(1,1,2)=d(2,1,1)=d(1,2,1)=16$, $d(1,2,2)=d(2,1,2)=d(2,2,1)=56$, $d(1,1,3)=d(3,1,1)=d(1,3,1)=24$,
$d(1,1,4)=d(4,1,1)=d(1,4,1)=32$, $d(1,2,3)=d(3,1,2)=d(2,3,1)=d(3,2,1)=d(1,3,2)=d(2,1,3)=128$, $d(2,2,2)=504$ for $r=3$. 
\\
\\
\noindent{\bf Remark.}
In spite of the negative terms 
in the formulas for the dimensions they give positive results. To understand why,
consider, for example, the case 
\begin{equation*}
d(2,2) = {\cal W}(2,2) - \frac{1}{2}{\cal W}(1,1)^2
\end{equation*}
with ${\cal W}(a,b)={\cal F}'=(a+b){\cal F}$, hence, $d(2,2)$ is four times the result in example 6. In the set of words counted by ${\cal F}'(2,2)=48$ there is a subset whose elements are words that are obtained gluing together  the words in the set counted by ${\cal W}(1,1)=8$. The gluing produces an overcounting which is corrected by the one half factor. So, 
$d(2,2)$ is positive. The same argument can be used to get positivity for the other formulas.

\begin{thm} \label{thm4.3}
For each $G_{r} \subseteq G_{R}$,
\begin{equation} \label{4.12}
\prod_{m_{i_{1}},...,m_{i_{r}}=1}^{\infty} 
(1+ z_{i_{1}}^{m_{i_{1}}}...z_{i_{r}}^{m_{i_{r}}}   )^{\theta_{+}} = 
e^{-g(z_{i_{1}}^{2}, ..., z_{i_{r}}^{2}) + 
g(z_{i_{1}}, ..., z_{i_{r}})}
\end{equation}
\begin{equation} \label{4.13}
\prod_{ m_{i_{1}},...,m_{i_{r}} =1 }^{\infty} 
(1-z_{i_{1}}^{m_{i_{1}}}...z_{i_{r}}^{m_{i_{r}}})^{\theta_{-}} =
e^{+g(z_{i_{1}}^{2}, ..., z_{i_{r}}^{2}) - g(z_{i_{1}}, ..., z_{i_{r}})}
\end{equation}
\end{thm}
\noindent {\bf Proof:} To prove \eqref{4.12} multiply and divide its left hand side by 
\begin{equation*}
\prod_{m_{i_{1}},...,m_{i_{r}}=1}^{\infty} 
(1- z_{i_{1}}^{m_{i_{1}}}...z_{i_{r}}^{m_{i_{r}}}   )^{\theta_{+}} 
\end{equation*}
 and use \eqref{4.8}.  
To get \eqref{4.13} write
\begin{equation*}
\prod_{ m_{i_{1}},...,m_{i_{r}} =1 }^{\infty} 
(1-z_{i_{1}}^{m_{i_{1}}}...z_{i_{r}}^{m_{i_{r}}})^{\theta_{-}} =
\prod_{N=r}^{\infty} 
\prod_{\begin{array}{c}
      \footnotesize{m_{i} >0}\\
       \footnotesize{m_{i_{1}}+...+m_{i_{r}}=N}
       \end{array}}
(1 - z_{i_{1}}^{m_{i_{1}}}...z_{i_{r}}^{m_{i_{r}}})^{\theta_{-}}
\end{equation*}
Decompose the product over $N$ into three products, namely,
one over all $N<2r$, one 
over all even $N \geq 2r$ and another one over all  odd $N > 2r$. Then, 
apply Theorem 2  and Theorem 3, formula \eqref{e17}. $\Box$

\begin{thm}\label{thm4.4}
\begin{equation} \label{4.14}
\prod_{ m_{i_{1}},...,m_{i_{r}}
=1 }^{\infty} 
(1+ z_{i_{1}}^{m_{i_{1}}}...z_{i_{r}}^{m_{i_{r}}} )^{\theta_{+}}
(1-z_{i_{1}}^{m_{i_{1}}}...z_{i_{r}}^{m_{i_{r}}} )^{\theta_{-}}
=1
\end{equation}
\end{thm}
\noindent {\bf Proof:} Multiply \eqref{4.12} and \eqref{4.13}. $\Box$
\\
\\
The left hand side of \eqref{A}  equals
\begin{equation*}
\prod_{j=1}^{R} (1 + z_{j})^{2} \prod_{r=2}^{R}\prod_{G_{r}}
\prod_{m_{i_{1}},...,m_{i_{r}} > 0}
(1 + z_{i_{1}}^{m_{i_{1}}} ...z_{i_{r}}^{m_{i_{r}}} )^{\theta_{+}}
(1-  z_{i_{1}}^{m_{i_{1}}} ...z_{i_{r}}^{m_{i_{r}}} )^{\theta_{-}}
\end{equation*}
From \eqref{4.14} Sherman identity follows.

\section{Two Lemmas} \label{sec : Lem}

\begin{lem}\label{l5.1} 
If
\begin{equation} \label{5.1}
g(n_{1},...,n_{k}) = \sum_{d \mid  n_{1},...,n_{k}} \frac{\mu(d)}{d}
f \left( \frac{n_{1}}{d},...,\frac{n_{k}}{d} \right)
\end{equation}
then
\begin{equation} \label{5.2}
f(n_{1},...,n_{k}) = \sum_{d \mid  n_{1},...,n_{k}} \frac{1}{d}
g \left( \frac{n_{1}}{d},...,\frac{n_{1}}{d} \right)
\end{equation}
\end{lem}
\noindent{\bf Proof:} Set
\begin{equation*}
G(n_{1},...,n_{k}) := (n_{1}+ ...+ n_{k}) g(n_{1},...,n_{k})
\end{equation*}
and
\begin{equation*}
F \left( \frac{n_{1}}{d},...,\frac{n_{k}}{d} \right)
:=
\left( \frac{n_{1}}{d}+...+\frac{n_{k}}{d} \right)
f \left( \frac{n_{1}}{d},...,\frac{n_{k}}{d} \right)
\end{equation*}
Then \eqref{5.1} can be expressed in the form
\begin{equation*}
G(n_{1},...,n_{k}) = \sum_{d \mid  n_{1},...,n_{k}} \mu(d)
F \left( \frac{n_{1}}{d},...,\frac{n_{1}}{d} \right)
\end{equation*}
M\"obius inversion gives
\begin{equation*}
F(n_{1},...,n_{k})= \sum_{d \mid  n_{1},...,n_{k}}
G \left( \frac{n_{1}}{d},...,\frac{n_{k}}{d} \right)
\end{equation*}
Therefore,
\begin{equation*}
(n_{1}+...+n_{k})f(n_{1},...,n_{k})
= \sum_{d \mid  n_{1},...,n_{k}}
\left( \frac{n_{1}}{d}+...+\frac{n_{k}}{d} \right) 
g \left( \frac{n_{1}}{d},...,\frac{n_{k}}{d} \right)
\end{equation*}
and the result follows.
The converse is also true.

\begin{lem}\label{l5.2}
 Let $N \geq \alpha =
n_{1}+...+n_{l}$, $n_{1},...,n_{l}$, $n_{i}>0$, a partition of
$\alpha$. Then,
\begin{equation} \label{5.3}
\sum_{\begin{array}{c}
       \footnotesize{\sum_{i=1}^{l} k_{i}=N}
       \end{array} } \prod_{i=1}^{l} \left( \begin{array}{c}
                   k_{i}-1 \\ n_{i}-1
\end{array}\right) =
\left( \begin{array}{c}
                   N-1 \\ \alpha-1
\end{array}\right)
\end{equation}
with the convention that a bracket in the left side is zero
whenever $k_{i}<n_{i}$.
\end{lem}

\noindent{\bf Proof:} Using
\begin{equation*}
\frac{q^{\alpha}}{(1-q)^{\alpha}} = \sum_{N=\alpha}^{\infty}
\left( \begin{array}{c}
                   N-1 \\ \alpha-1
\end{array}\right)
q^{N}
\end{equation*}
it follows that
\begin{eqnarray*}
\frac{q^{\alpha}}{(1-q)^{\alpha}} &=& \prod_{i=1}^{l}
\frac{q^{n_{i}}}{(1-q)^{n_{i}}}
= \prod_{i=1}^{l}
\sum_{k_{i}=n_{i}}^{\infty} \left(
\begin{array}{c}
                   k_{i}-1 \\ n_{i}-1
\end{array}\right)
q^{k_{i}}\\
&=& \sum_{N = \alpha}^{\infty} \sum_{\begin{array}{c}
      \footnotesize{k_{i} \geq n_{i}}\\
       \footnotesize{\sum_{i=1}^{l} k_{i}=N}
       \end{array} } \prod_{i=1}^{l} \left(
\begin{array}{c}
                   k_{i}-1 \\ n_{i}-1
\end{array}\right)
q^{N}
\end{eqnarray*}
Comparison with previous expression and the convention gives the
result.

\noindent{\bf Acknowledgements.} We kindly thank  Prof. Peter Moree (Max Planck Institute for Mathematics, Bonn) and Prof. Thomas Ward (University of East Anglia, UK) for email correspondence regarding the positivity of M\"obius inversion formula.

\end{document}